\documentclass{elsart} 

\usepackage{t1enc}
\usepackage[final]{graphicx}
\usepackage{graphicx}
\usepackage{epsfig}
\usepackage{amssymb}

\usepackage{color}

\begin{document}

\newcommand{\ii}{\'{\i}}
\newcommand{\be}{\begin{equation}}
\newcommand{\ee}{\end{equation}}
\newcommand{\ben}{\begin{eqnarray}}
\newcommand{\een}{\end{eqnarray}}
\newcommand{\pp}{\prime}
\newcommand{\nn}{\nonumber}
\newcommand{\nd}{\noindent}
\newcommand{\tr}{\textcolor{red}}

\begin{frontmatter}

\title{A Shannon-Tsallis transformation}

\author{E. Rufeil Fiori$^1$   }
\address{Facultad de Matem\'{a}tica, Astronom\'{i}a y F\'{i}sica, Universidad Nacional de C\'{o}rdoba \\
                Ciudad Universitaria, 5000 C\'{o}rdoba, Argentina}

\author{A. Plastino$^{2,\,3}$  }
\ead{angeloplastino@gmail.com}
\address{ $^{2}$National University La Plata
\& CONICET IFLP-CCT, C.C. 727 - 1900 La Plata, Argentina
\\ $^{3}$ Universitat de les Illes Balears and IFISC-CSIC, 07122 Palma de Mallorca, Spain}

\begin{abstract}
\noindent We determine a general link between two different solutions of
 the MaxEnt variational problem, namely, the ones that correspond to using either Shannon's
 or Tsallis' entropies in the concomitant variational problem. It is
 shown that the two variations lead to equivalent solutions that  take different appearances but contain the same information.
 These solutions are linked by our transformation.
\end{abstract}


\end{frontmatter}

\section{Introduction}

\nd Nonextensive statistical mechanics (NEXT) \cite{[1],[2],AP}, a
 generalization of the orthodox Boltzmann-Gibbs (BG) one, is actively
investigated and applied in many areas of scientific endeavor.
NEXT is based on a nonadditive (though extensive \cite{[3]})
entropic information measure, that is characterized by a real
index q (with q = 1 recovering the standard BG entropy). It has
been used with regards to  variegated systems such as cold atoms
in dissipative optical lattices \cite{[4]}, dusty plasmas
\cite{[5]}, trapped ions \cite{[6]}, spinglasses \cite{[7]},
turbulence in the heliosheath \cite{[8]}, self-organized
criticality \cite{[9]}, high-energy experiments at LHC/CMS/CERN
\cite{[10]} and RHIC/PHENIX/Brookhaven \cite{[11]},
low-dimensional dissipative maps \cite{[12]}, finance \cite{[13]},
galaxies \cite{AP1}, Fokker-Planck equation's applications
\cite{AP2}, etc.

\nd A typical NEXT feature is that it can be  expressed by
recourse to generalizations \`a la q of standard mathematical
concepts \cite{borges}. Included are, for instance, the logarithm
[q-logarithm] and exponential functions (usually denoted as
$e_q(x)$, with $e_{q=1}(x)=e^{x}$), addition and multiplication,
Fourier transform (FT) and the Central Limit Theorem (CLT)
\cite{tq2}. The q-Fourier transform $F_q$ exhibits the nice
property of
 transforming q-Gaussians into q-Gaussians \cite{tq2}.
Recently, plane waves, and the representation of the Dirac delta
into plane waves have been also generalized
\cite{[15],[16],tq1,tq4}.

\nd Our central interest here resides in the q-exponential
function, regarded as the MaxEnt variational solution \cite{AP} if
the pertinent  information measure is Tsallis' one. We will show
that there is a transform procedure that converts any
Shannon-MaxEnt solution \cite{jaynes,katz} into a q-exponential,
without modification of the associated Lagrange multipliers, that
carry with them all the physics of the problem at hand. Why?
Because of the Legendre transform properties of the MaxEnt
solutions [see for instance
\cite{jaynes,katz,mine1,mine2,mine3,mine4,mine5,pp97} and
references therein].

\nd Accordingly, we are here proving that the physics of a given
problem can be discussed in equivalent fashion  by recourse to
either Shannon's measure or Tsallis' one, indistinctly.

\section{The central idea}

\nd We wish to connect orthodox exponentials with q-exponentials.
Let us consider the Shannon-MaxEnt solution for a constraint given
by the average value of the variable in question, that we call $u$:

\ben \label{uno}
    & p(u)du= \exp{[-\mu-\lambda u]}\,du=e^{-\mu}\,e^{-\lambda u}du, \cr
    & \langle u \rangle= K, \cr
    & \int p(u)du=1,
\een where the two Lagrange multipliers $\mu$ and $\lambda$
correspond, respectively, to normalization and conservation of the
$u-$mean value $<u>=K$.

\nd Consider now a second variable $x$ such that $dx/du=g(x)$,
with $g(x)$ a function we will wish to determine below.  Assume
that in the second variable we express the Tsallis-MaxEnt
solution, with the same constraints, but employing the above
mentioned q-exponential functions ($e_q(x)= [1+(1-q)x]^{1/(1-q)}$.
The support of this function is sometimes finite, depending on the
$q-$value. See more details in, for instance, \cite{[2]}).

\be \label{dos}
    p(x)dx= C\,e_q(-\lambda\,x)\,dx;\,\,\,C= {\rm normalization \,\,const.}
\ee

\nd We want $$p(x)dx=p(u)du.$$

\nd This entails:

\be \label{tres}
    C\,e_q(-\lambda\,x)\,dx = p(x)dx= e^{-\mu}\,e^{-\lambda u(x)} (dx/g(x)),
\ee
that is, given that $du= g(x)^{-1}\,dx$,

\be \label{cuatro}
    C\,g(x)\, e_q(-\lambda\,x)=e^{-\mu}\,\exp{[-\lambda \int\,dx \,g(x)^{-1}]}.
\ee
Remember now that

\be \label{cinco}
    \frac{d\,e_q(x)}{dx}= e_q(x)^{q},
\ee
so that, taking the logarithm of Eq. (\ref{cuatro}) we find

\be \label{seis}
    \ln{C} + \ln{g(x)} + \ln{[e_q(-\lambda\,x)]}= -\mu - \lambda \, \int dx \,g(x)^{-1}.
\ee

\nd Now, we derive w.r.t. $x$ Eq. (\ref{seis}) and have

\be \label{siete}
    \frac{g'(x)}{g(x)} -\lambda\,e_q(-\lambda\,x)^{q-1} = -\lambda \,g(x)^{-1},
\ee
which leads to a differential equation for our desired $g(x)$:

\be \label{ocho}
    g'(x) -\lambda\, e_q(-\lambda\,x)^{q-1}\,g(x)+\lambda=0.
\ee
\nd
Solving this equation we establish the link we are looking for.

\section{The differential equation}

\nd For simplicity we now set

\be \label{nueve}
    P(x)= -\lambda\, e_q(-\lambda\,x)^{q-1};\,\,\,Q(x)=-\lambda,
\ee
and recast (\ref{ocho}) as

\be \label{diez}
    g'(x) + P(x)\,g(x) = Q(x).
\ee
Introduce now the integrating factor $I$

\be \label{once}
    I= \exp{[\int\,dt\,P(t)]},
\ee
and multiply (\ref{diez}) by it

\be \label{doce}
    g'(x)I + P(x)I\,g(x) = Q(x)I.
\ee
Note that

\be \label{trece}
    \frac{d}{dx}[g(x)I]= g'(x)I+  P(x)I\,g(x),
\ee
so that (\ref{once}) becomes

\be \label{catorce}
    \frac{d}{dx}[g(x)I]= Q(x)I.
\ee
Integrating this we have now

\be \label{quince}
    g(x)I= Q(x) \int Idx + c.
\ee
Finally, we can formally express our ``solution" function $g(x)$ as

\be \label{dieciseis}
    g(x)= \frac{ \int Q(x)Idx + c}{I}.
\ee
\nd
Thus, for the differential equation

\be \label{diecisiete}
    g'(x)-\lambda e_q(-\lambda x)^{(q-1)}g(x)+\lambda = 0,
\ee
we obtain the solution

\ben \label{dieciocho}
    g(x)&=& \exp\left(\lambda \int^x e_q(-\lambda x')^{(q-1)}\mathrm{d}x' \right), \nonumber \\
    &&\left[-\lambda \int^x \exp \left( -\lambda \int^{x'} e_q(-\lambda x'')^{(q-1)}\mathrm{d}x''\right)\mathrm{d}x'+c\,\right],
\een
where c is an integration  constant.  Now, using

\be \label{diieciocho}
   \lambda \int^x e_q(-\lambda x')^{(q-1)}\mathrm{d}x' = \ln\left(e_q(-\lambda x)^{(-1)}\right),
\ee
and

\be \label{diiecinueve}
       -\lambda \int^x \exp \left( -\lambda \int^{x'} e_q(-\lambda x'')^{(q-1)}\mathrm{d}x''\right)\mathrm{d}x' = \frac{e_q(-\lambda x)^{(2-q)}}{2-q},
\ee
we finally arrive at

\ben \label{veinte}
    g(x) &=& e_q(-\lambda x)^{(-1)} \left[ \frac{e_q(-\lambda x)^{(2-q)}}{2-q} + c \right],
\een
which provides us with the Tsallis-Shannon Jacobian

\be \label{veintiuno}
    J(x)= \frac{1}{g(x)}.
\ee

Consider now the instance $q\rightarrow1$. One obviously ought to have $g(x)=1$. We face
\ben \label{veintiuno1}
    g(x)&\rightarrow& e^{\lambda x} \left[ e^{-\lambda x} + c \right] \\
       &=& 1+c\,e^{\lambda x},
\een
which entails $c=0$ and

\ben \label{veintiuno2}
    g(x) &=& \frac{e_q(-\lambda x)^{(1-q)}}{2-q} = \frac{1-(1-q)\lambda\,x}{2-q}.
\een

\subsection{Expansion near q=1}

\nd Let us now take $q=1-\epsilon$ in the solution $g(x)=g(x;q)$ (Eq.(\ref{veintiuno2})),

\be \label{veintidos}
    g(x;q=1-\epsilon)= \frac{1-\epsilon \lambda x}{1+\epsilon}.
\ee
Thus, a first-order  expansion in $\epsilon$ gives:

\ben \label{veintitres}
    g(x;q=1-\epsilon)&=& 1-\left(1+\lambda x \right) \epsilon +O(\epsilon^{2}).
\een

\subsection{The case $q=2$}

\nd Now, let us look at the $q = 2-\epsilon$ scenario for our solution $g(x)=g(x;q)$;

\be \label{veinticuatro}
    g(x; q = 2 -\epsilon) = \frac{1+(1-\epsilon)\lambda x}{\epsilon}.
\ee
That, in the limit of $\epsilon\rightarrow 0$ we have

\ben
    g(x;q=2-\epsilon) &\asymp& \frac{(1+\lambda x)}{\epsilon},\\
                    &\rightarrow& +\infty \,\,\, \mathrm{as} \,\,\, \epsilon\rightarrow 0^{+} \,\,\, (\mathrm{i.e.,}\,q\rightarrow 2^{-})\,\,\,\mathrm{and}\\
                    &\rightarrow& -\infty \,\,\, \mathrm{as} \,\,\, \epsilon\rightarrow 0^{-} \,\,\, (\mathrm{i.e.,}\,q\rightarrow 2^{+}).
\een
There is a divergence in the first term, in the form of $1/\epsilon$.

\nd We are then in a position to state that, for $\lambda x > -1$,
and given the definition of the q-exponential (see \cite{[2]}), our
``inverse-Jacobian" function $g(x)$ is positive for $q<2$ and
negative for $q>2$, while diverging at $q=2$. We conclude that the
transform we are studying is not valid only in the isolated case
$q=2$ and changes sign there.

\section{Arbitrary constraint}

\nd We generalize now the preceding considerations to the case of
a generalized constraint $<~h(x)>$, with $h \in \mathcal{L}_2$.
The concomitant  Shannon MaxEnt solution is \cite{katz}

\ben  \label{veintiseis}
    p(u)du &=& e^{-\mu}\,e^{-\lambda h(u)}du, \\
    \langle h(u) \rangle &=& K, \\
    \int p(u)du &=& 1.
\een
while the  Tsallis-MaxEnt solution with the same constraint becomes

\ben \label{veintisiete}
    p(x)dx &=& C\,e_q(-\lambda\,h(x))\,dx;\,\,\,C= {\rm normalization \,\,const.}
\een
Assume that $u(x)$ exists and call, as before, $dx/du=g(x)$. We have, as our cornerstone the relation

\be \label{veintiocho}
    p(x)dx=p(u)du,
\ee
entailing

\be   \label{veintinueve}
    C\,e_q(-\lambda\,h(x))\,dx = e^{-\mu}\,e^{-\lambda h(u(x))} (dx/g(x)),
\ee
so that, taking the logarithm to this equation we find

\be \label{treinta}
    \ln{C} + \ln{[g(x)]} + \ln{[e_q(-\lambda\,h(x))]}= -\mu - \lambda \, h(u(x)).
\ee
Taking derivatives   w.r.t. $x$ yields

\be \label{treintiuno}
    \frac{g'(x)}{g(x)} - \lambda\,e_q(-\lambda\, h(x))^{q-1} h'(x)= -\lambda \,\frac{1}{g(x)}h'(x),
\ee
which leads to a differential equation for our desired transformation function $g(x)$

\be \label{treintidos}
    g'(x) -\lambda\, e_q(-\lambda\,h(x))^{q-1}\,h'(x)g(x)+\lambda h'(x)=0,
\ee
quite similar in shape to Eq. (\ref{diecisiete}), being thus solved in similar fashion. Using  Eq. (\ref{dieciseis}) we
encounter

 \ben \label{treintitres}
    \int P &=& \int e_q(-\lambda h(x))^{q-1} (-\lambda h'(x))\mathrm{d}x, \\
    &=& \ln(e_q(-\lambda h(x))), \\
    e^{-\int P}&=& e_q(-\lambda h(x))^{-1}, \\
    \int Q e^{\int P} &=& \int e_q(-\lambda h(x)) (-\lambda h'(x))\mathrm{d}x, \\
     &=& \frac{e_q(-\lambda h(x))^{2-q}}{2-q},
\een
that leads to

\be \label{treinticuatro}
    g(x)=e_q(-\lambda h(x))^{-1} \left[ \frac{e_q(-\lambda h(x))^{2-q}}{2-q} +
    c\right].
\ee It is clear that for $q=1$ we have $g(x)=1$ if c=0. Thus,

\be \label{treinticuatro1}
    g(x)=\frac{1-(1-q)\lambda h(x)}{2-q}.
\ee
Also, the regime-change at $q=2$ discussed above does not change.

\subsection{Special case $h(x)=x^{2}$}

\nd For the special case of a variance constraint we have

\be  \label{treinticinco}
    g'(x) -\lambda\, e_q(-\lambda\,x^{2})^{q-1}\,2xg(x)=-\lambda 2x,
\ee
whose solution is

\be \label{treintiseis}
    g(x)=\frac{1-(1-q)\lambda x^{2}}{2-q}.
\ee
Let us take now  $q=1-\epsilon$ in the solution $g(x)=g(x;q)$. Then, an expansion near $\epsilon=0$ yields

\ben
    g(x;q=1-\epsilon)&=& 1 - \left(1+\lambda x^{2}\right) \epsilon +O(\epsilon^{2}).
\een

\section{Generalization to M constraints}

Now, we generalize to the case of $M$ constraints of the form:
\be
     \langle h_i(u) \rangle = K_i;\,\,\,\, i=1,\ldots,M.
\ee
Lets use vector notation and call $\mathbf{h}=(h_1,h_2,\ldots,h_M)$ and $\mathbf{\lambda}=(\lambda_1,\lambda_2,\ldots,\lambda_M)$ its associated Lagrange multipliers. Then we have for the Shannon MaxEnt solution

\be
    p(u)du = e^{-\mu}\,e^{-\mathbf{\lambda}\cdot\mathbf{h}(u)}du,
\ee
while the  Tsallis-MaxEnt solution with the same constraints becomes

\ben
    p(x)dx &=& C\,e_q(-\mathbf{\lambda}\cdot\mathbf{h}(x))\,dx;\,\,\,C= {\rm normalization \,\,const.}
\een
Assume that $u(x)$ exists and call, as before, $dx/du=g(x)$. Then, following the steps of the previous section and solving the  corresponding differential equation, we obtain

\be
    g(x)=e_q(-\mathbf{\lambda}\cdot\mathbf{h}(x))^{-1} \left[ \frac{e_q(-\mathbf{\lambda}\cdot\mathbf{h}(x))^{2-q}}{2-q} +
    c\right],
\ee where $c=0$ in order to get $g(x)=1$ in the limit
$q\rightarrow 1$. Thus,

\be
    g(x)=\frac{1-(1-q)\mathbf{\lambda}\cdot\mathbf{h}(x)}{2-q}.
\ee

\section{Conclusions}

\nd We have shown here  that, from a MaxEnt practitioner
view-point, one can indistinctly employ Shannon's logarithmic
entropy $S$ or Tsallis' power-law one $S_q$ (for any $q$ except
$q=2$). The physics described is the same. To choose between $S$
and $S_q$ is just a matter of convenience in the sense of getting
simpler expressions in one case than in the other.

\nd The link between the two concomitant probability distributions
$P_{Shannon}(x)$ and  $P_{Tsallis}(x)$   is given by the Jacobian
$J=1/g$, where $g$ is the simple function

\be
    g(x)=\frac{1-(1-q)\lambda h(x)}{2-q},
\ee
with $\lambda$ the pertinent Lagrange multiplier, and $h(x) \in \mathcal{L}_2$ an arbitrary
function whose mean value $<h>$ constitutes MaxEnt's informational
input.

\section{Appendix: The four different Tsallis' treatments}

\nd A savvy Tsallis practitioner may wonder what happens with the
four different ways of computing $q-$mean values that one finds in
Tsallis' literature (see \cite{ferri4} and references therein). In
addition to the normal expectation values we have employed above,
one also encounters, for a quantity $A(x)$ averaging ways that,
themselves, depend upon $q$ (see below). This transforms our
cornerstone-equality $p_{Shannon}(u)du=p_{Tsallis}(x)dx$ into
something much more complicated. However, there is a way out,
following the discoveries reported in \cite{ferri4}.

\nd  Bernoulli published in the {\em Ars Conjectandi} the first
formal attempt to deal with probabilities already  in 1713 and
Laplace further formalized the subject in his {\em Th\'eorie
analytique des Probabilit\'es} of 1820. In the intervening
centuries Probability Theory (PT) has grown into a rich, powerful,
and extremely useful branch of Mathematics. Contemporary Physics
heavily relies on PT for a large part of its basic structure,
Statistical Mechanics \cite{patria,reif,sakurai}, of course, being
a most conspicuous example. One of PT  basic definitions is that
of the mean value of an observable $\mathcal{A}$ (a measurable
quantity). Let $A$ stand for the linear operator or dynamical
variable associated with $\mathcal{A}$. Then,

\be
    \left<A\right>=\int\,dx\,p(x)\,A(x).
\ee

This was the averging procedure that Tsallis used in his first, pioneering 1988 paper
\cite{[1]}, and the one discussed in the preceding Sections. It is
well known that, in some specific cases, it  becomes necessary to
use  ``weighted" mean values, of the form
\be
    \left<A\right>=\int\,dx\,f[p(x)]\, A(x),
\ee
with $f$ an analytical function of  $p$. This happens, for instance, when
  there is a  set of states characterized by a distribution  with a
  recognizable maximum and a large tail that contains low but finite probabilities.
One faces then the need of making  a  pragmatical (usually of
experimental origin) decision regarding $f$ \cite{ferri4}.
\nd In the first stage of NEXT-development, its pioneer practitioners made the
pragmatic choice of using ``weighted" mean values, of rather
unfamiliar appearance for many physicists. Why? The reasons were
of theoretical origin. It was at the time believed that, using the
familiar linear, unbiased mean values, one was unable to get rid
of the Lagrange multiplier associated to
probability-normalization. Since the Tsallis' formalism yields, in
the limit $q\rightarrow 1$, the orthodox Jaynes-Shannon
 treatment, the natural  choice was to construct weighted expectation values (EVs) using the
  index  $q$,

\be
    \left<A\right>_q=\int\,dx\,p(x)^q A(x),
\ee
the so-called Curado-Tsallis unbiased mean values (MV) \cite{CT}. As shown in
\cite{AP}, employing the Curado-Tsallis (CT) mean values allowed
one to obtain an {\it analytical} expression for the partition
function out of the concomitant MaxEnt process \cite{AP}. This EV
choice leads to a non extensive formalism endowed with interesting
 features: i) the above mentioned property of its partition
 function $Z$, ii)
a  numerical treatment that is relatively simple, and iii) proper
results in the limit $q\rightarrow 1$.  It has, unfortunately, the
 drawback of exhibiting  un-normalized mean values, i.e.,
$\left<\left<1\right>\right>_q\neq 1$. The latter problem was
circumvented in the subsequent work of Tsallis-Mendes-Plastino
(TMP) \cite{TMP}, that ``normalized" the CT treatment by employing
mean values of the form

\be
    \left<A\right>_q=\int\,dx \,\frac{p(x)^q}{{\mathcal X}_q}\,
A(x);\,\,\,{\mathcal X}_q= \int\,dx \,p(x)^q.
\ee
 Most NEXT works employ the TMP procedure. However, the concomitant
treatment is not at all simple. Numerical  complications often
ensue, which has encouraged the development of different,
alternative  approach called the OLM one \cite{OLM}, that
preserves the main TMP-idea  (the ${\mathcal X}_q$ normalization
sum) but is numerically simpler. Now, despite appearances, the
four Tsallis' treatments are equivalent, as shown in
\cite{ferri4}. By equivalence we mean that if one knows the
probability treatment $P_i; \,\,i=1,2,3,4$ obtained by anyone of
the four treatments, there is a unique, automatic way to write
down $P_j; \,\,j\ne i$. More precisely, any $P_i$ is a
q-exponential, and they all possess the same information-amount if
the pair $q,\,\beta$ is appropriately ``translated'' from a version
to the other \cite{ferri4}. Indeed,

\ben
    &  P_{i}=Z^{-1}\exp _{q^{\ast }}(-\beta ^{\ast }x),\cr
    & Z=\sum_{i}\exp _{q^{\ast }}(-\beta ^{\ast }x).
\een
Then, as made explicit in  \cite{ferri4}, given any of the four possible $P_i$'s
$q^{\ast},\,\,\beta ^{\ast}$, one can get the $q,\,\,\beta$ values
appropriate for the q-exponential of any $j \ne i$. Such
``dictionary'' allows one to translate the results obtained in the
preceding Sections to any other of the three remaining averaging
procedures.

\vskip 6mm

\nd {\bf Acknowledgments}  This work was partially supported by:
i) the project PIP1177 of CONICET (Argentina) and  ii)  the
project FIS2008-00781/FIS (MICINN) and FEDER (EU) (Spain, EU).

\end{document}